\documentclass[12pt]{article}
\input epsf
\usepackage[T2A,OT1]{fontenc}
\usepackage[ot2enc]{inputenc}
\usepackage[russian,english]{babel}
\usepackage{epsfig}
\usepackage{psfig}
\input cyracc.def

\def\L{{\mathcal{L}}}
\def\Z{{\mathbf{Z}}}
\def\R{{\mathbf{R}}}
\def\sgn{{\mathrm{sgn}\,}}
\def\const{{\mathrm{const}}}

\def\ch{{\mathrm{ch\,}}}

\title{Uniform approximation of $\sgn x$ by polynomials and entire functions}
\author{Alexandre Eremenko\thanks{Supported by NSF grants DMS-0555279, DMS-244547.}\,
and Peter
Yuditskii\thanks{Partially supported by {\em Marie Curie Intl. Fellowship}
within the 6-th EC Framework Progr. Contract MIF1-CT-2005-006966.}}
\begin{document}
\maketitle

\vspace{.2in}

In 1877, E. I. Zolotarev \cite{Zol,Akh}
found an explicit expression, in terms of elliptic functions,
of the rational function
of given degree $m$ which is uniformly closest to $\sgn(x)$
on the union of two intervals $[-1,-a]\cup[a,1]$.
This result was subject to many generalizations, and it has
applications in electric engineering. 

Surprisingly, to the best of our knowledge,
the similar problem for polynomials
was not solved yet, so we investigate it in this paper.

For comparison, we mention here the results on the
uniform approximation of $|x|^\alpha,\;\alpha>0$ on $[-1,1]$.
Polynomial approximation
was studied by S.~Bernstein \cite{Bern1,Bern2}
who found that for the error $E_m(\alpha)$ of the best approximation
by polynomials of degree $m$ the following limit exists:
$$\lim_{m\to\infty}m^\alpha E_m(\alpha)=\mu(\alpha)>0.$$
This result for $\alpha=1$ was obtained by Bernstein in 1914,
and he asked the question, whether one can express
$\mu(1)$ in terms of some known transcendental functions.
This question is still open. Bernstein also obtained in \cite{Bern2}
the asymptotic relation $$\lim_{\alpha\to 0}
\mu(\alpha)=1/2.$$

The analogous problem of uniform rational approximation of $|x|^\alpha,\;\alpha>0$
on $[-1,1]$ 
was recently solved by H. Stahl \cite{Sta}, who completed a long
line of development with a remarkably
explicit answer:
$$\lim_{m\to\infty}\exp(\pi\sqrt{\alpha m})E^r_m=
2^{2+\alpha}|\sin(\pi\alpha/2)|,$$
where $E^r_m$ is the error of the best rational approximation. 

Now we state our results. 
Let $p_m$ be the polynomial
of degree at most $2m+1$
of least deviation from $\sgn(x)$ on $X(a)=[-1,-a]\cup[a,1]$ where
$0<a<1$. It follows from the general theory of Chebyshev that such
polynomial  is unique. Put $L_m(a)=\max_{X(a)}|p_m(x)-\sgn(x)|.$ Then we have 
\vspace{.1in}

\noindent
{\bf Theorem 1} {\em The following limit exists
$$\lim_{m\to\infty}
\sqrt{m}\left(\frac{1+a}{1-a}\right)^{m}L_m(a)=\frac{1-a}{\sqrt{\pi a}}.$$
}
\vspace{.1in}

{\em Remark.} When approximating an odd function on a symmetric
set, polynomials of even degrees are useless. Indeed, if $q$ is
a polynomial of degree at most $2m$ which deviates least from our function,
then
$(q(z)-q(-z))/2$ is an odd polynomial, thus its degree is less
than $2m$ and its deviation is at most that of $q$. Thus $q$ is of odd
degree.
\vspace{.1in}

Our approximation problem is equivalent to a problem of
weighted approximation on a single interval. Indeed,
$p_m$ can be written as $p_m(x)=xq_m(x^2)$, where $q_m$
is the polynomial of degree $m$ that minimizes the weighted uniform
distance
\begin{equation}
\label{our}
\sup_{[a^2,1]}\sqrt{x}|q(x)-1/\sqrt{x}|.
\end{equation}
over all polynomials $q$ of degree at most $m$.

It is useful to compare our result with the result of Bernstein
\cite{Bernstb}, see also \cite[Additions and Problems, 44]{Akha}
that gives the rate of the best {\em unweighted} uniform polynomial
approximation of $1/\sqrt{x}$:
\begin{equation}
\label{bern}
\lim_{m\to\infty}\sqrt{m}\left(\frac{1+a}{1-a}\right)^m\inf_{\deg r=m}
\,\sup_{[a^2,1]}|r(x)-1/\sqrt{x}|=\frac{1}{2\sqrt{\pi}}(1-a^2)a^{-3/2}.
\end{equation}

In our proof of Theorem 1, the asymptotics of the error term
is obtained in the form
\begin{equation}
\label{A}
\lim_{m\to\infty}\sqrt{m}\left(\frac{1+a}{1-a}\right)^mL_m(a)=
e^{-c}\frac{\sqrt{2}(1-a)}{\sqrt{a}},
\end{equation}
with
\begin{equation}
\label{int}
c=\frac{1}{\pi}\int_0^\infty\left(\Im H(t)-\frac{\pi}{2}
\chi_{[2,\infty)}\right)\frac{dt}{t},
\end{equation}
where $\chi_{[2,\infty)}$ is the characteristic function of
the ray $[2,\infty)$, and $H$ is the 
conformal map of the upper half-plane onto the region
in the upper half plane above the curve
$$\{ t+i\arccos e^{-t}:t\geq 0\},$$
normalized by $H(0)=0$ and $H(z)\sim z,\, z\to\infty$.
Then the numerical value  $c=(1/2)\log(2\pi)$ is derived
from comparison of (\ref{A}) with (\ref{bern}) for $a\to 1$.
So, as a curious corollary from Theorem 1 and the result
of Bernstein, we evaluate the integral (\ref{int}).

Following Bernstein, we also consider approximation by entire functions
of exponential type. Let $L(A)$ be the error of the best uniform approximation
of $\sgn(x)$ by entire functions of exponential type one, on the set
$(-\infty,-A]\cup [A,+\infty)$. 
\vspace{.1in}

\noindent
{\bf Theorem 2} {\em The following limit exists
$$\lim_{A\to\infty}
\sqrt{A}\exp(A)L(A)=\sqrt{2/\pi}.$$
}
\vspace{.1in}

The proof of Theorem 2 is similar to (and simpler than) that of Theorem~1.
On the best $L^1$ approximation of $\sgn(x)$ by entire functions
of exponential type we refer to \cite{Vaal}.

Our proofs are based on special representations of polynomials and entire
functions of best approximation which are of independent interest.
It follows from Chebyshev's theory (see, for example, \cite[Ch. II]{Akha})
that polynomials $p_m$ are characterized
by the property that the difference $p_m(x)-\sgn(x)$ takes its extreme
values $\pm L_m(a)$ on $X(a)$ $2m+4$ times intermittently, so that the
graph of $p_m$ looks like this:
\begin{center}
\epsfxsize=2.0in%
\centerline{\epsffile{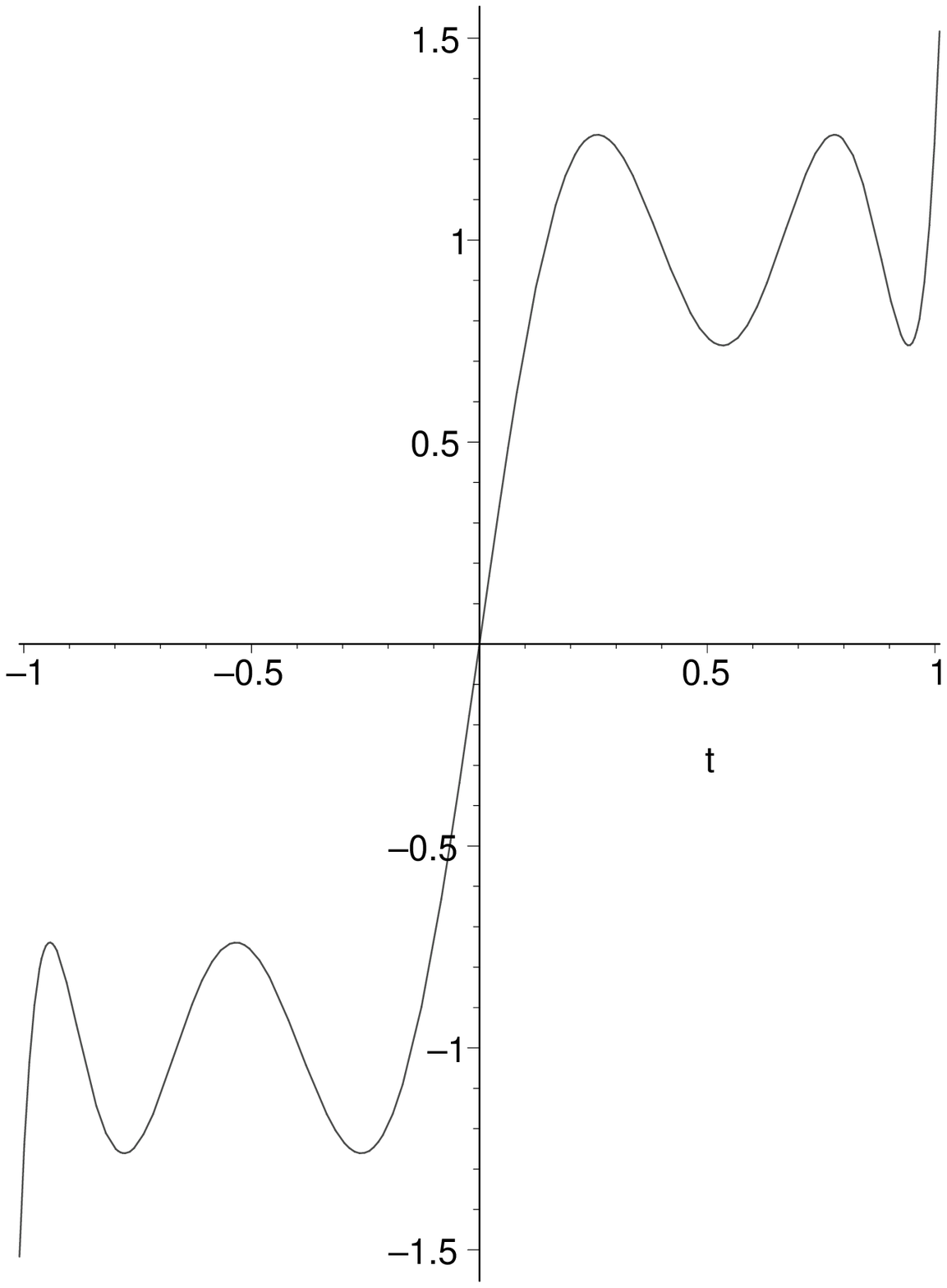}}
\nopagebreak
Fig. 1. Graph of $p_4$ with $a=0.1$.
\end{center}

The extremal entire function is unique and is characterized
by the properties that it has no asymptotic values,
all its critical
values are real; those on the negative ray are $-1\pm L(A)$
and those on the positive ray are $1\pm L(A)$.

We mention a general theorem of Maclane \cite{Mc} and Vinberg \cite{Vin}
on the existence and uniqueness of real polynomials and entire functions
with prescribed
(ordered!) sequences of critical values.
In the case of polynomials,
this theorem says that there is a one to one correspondence between
finite ``up-down'' real sequences $$\ldots \leq c_{k-1}\geq c_k\leq
c_{k+1}\geq
\ldots,$$ and real polynomials whose all
critical points are real, modulo a change of the independent
variable $z\mapsto az+b,\; a>0, b\in\R$.
If $(x_k)$ is the sequence of critical points of such polynomial, then
$c_k=P(x_k)$.

The MacLane--Vinberg theorem
is based on an explicit description of the Riemann surfaces
spread over the plane of the inverse functions $P^{-1}$. 
This explicit description and
our previous work \cite{Ere}, \cite{SY} suggest to look for a representation
of these extremal polynomials and entire functions
in the form $\cos\phi(z)$ where $\phi$ is an appropriate
conformal map.  

Let $L$ and $B$ be positive numbers, $0<L<1$, and
\begin{equation}
\label{BL}
L=\frac{1}{\ch B}\sim 2e^{-B}, \quad B\to\infty.
\end{equation}
Consider the component $\gamma_B\subset\{ z:\Re z\in [0,\pi],\,\Im z>0\}$
of the preimage of the ray
$$\{ w:\Re w=1/L,\; \Im w<0\}$$
under $w=\cos z$.
It is easy to see that this curve $\gamma_B$ can be parametrized
as 
$$\gamma_B= \{ \arccos(\ch B/\ch t)+it: B\leq t<\infty\}.$$
This curve begins at $iB$ and then goes to infinity approaching the
line $\{ \pi/2+it, t>0\}$ with exponential rate. 

Let $\Omega_B$ be the region in the upper half-plane
whose boundary consists of the positive ray, the vertical segment $[0,iB]$
and the curve $\gamma_B$. For fixed $B>0$, let $\phi_B$ be the conformal map
of the first quadrant onto $\Omega_B$ such that $\phi_B(z)\sim z,\; z\to\infty$
and $\phi_B(0)=iB$. Let $A=A(B)=\phi_B^{-1}(0)$. Then $A$ is a continuous 
strictly increasing function of $B$, and we may consider the inverse function
$B(A)$. 
\vspace{.1in}

\noindent
{\bf Theorem 3} {\em The approximation error in Theorem 2 is 
$L(A)=1/\ch B(A)$, and the extremal function
can be defined in the first quadrant
by the formula}
$$1-L(A)\cos\phi_{B(A)}.$$
\vspace{.1in}

Let $\Omega_{B,m}$ be the region in the
half-strip
$\{ z:\Re z\in (0,\pi(m+1)),\,\Im z>0\}$ bounded on the left by $\gamma_B$.
Let $\phi_{B,m}$ be the conformal map of the first quadrant onto $\Omega_{B,m}$
such that $\phi_{B,m}(0)=iB,\; \phi_{B,m}(1)=\pi(m+1)$ and $\phi_{B,m}(\infty)=\infty.$
Let $a=a(B,m)=\phi_{B,m}^{-1}(0).$ Then $a(B,m)$
is a continuous increasing function
of $B$ for fixed $m$,
so it has the inverse $B_m(a)$.

\vspace{.1in}

\noindent
{\bf Theorem 4} {\em The error term in Theorem 1 is 
$L_m(a)=1/\ch B_m(a),$ and the extremal polynomial is given in
the first quadrant by
$$1-L_m(a)\cos\phi_{B,m},\quad\mbox{where}\quad B=B_m(a).$$}

Discontinuous functions cannot be uniformly approximated by
polynomials with arbitrarily small precision,
however, in our situation we can obtain an approximation
which seems to be the second best
thing to the uniform approximation.

We introduce the notation\footnote{If $[-1,1]$ is replaced by $(-\infty,
\infty)$ this becomes the L\'evy distance. It is really a distance on the
set of bounded increasing functions on the real line \cite[Ch. VIII]{LO}.}
$$\L(f,g)=\inf\{ h:f(x-h)-h\leq g(x)\leq f(x+h)+h,\; -1\leq x\leq 1 \}.$$
Thus the statement that $\L(\sgn,g)\leq \epsilon$ means that the graph
of the restriction of $g$ on $[-1,1]$ belongs to a ``rectangular
corridor'' of width $\epsilon$ around the ``completed graph'' of $\sgn(x)$, which consists
of the graph of $\sgn(x)$ and the vertical segment $[(0,-1),(0,1)]$.
\begin{center}
\epsfxsize=3.0in%
\centerline{\epsffile{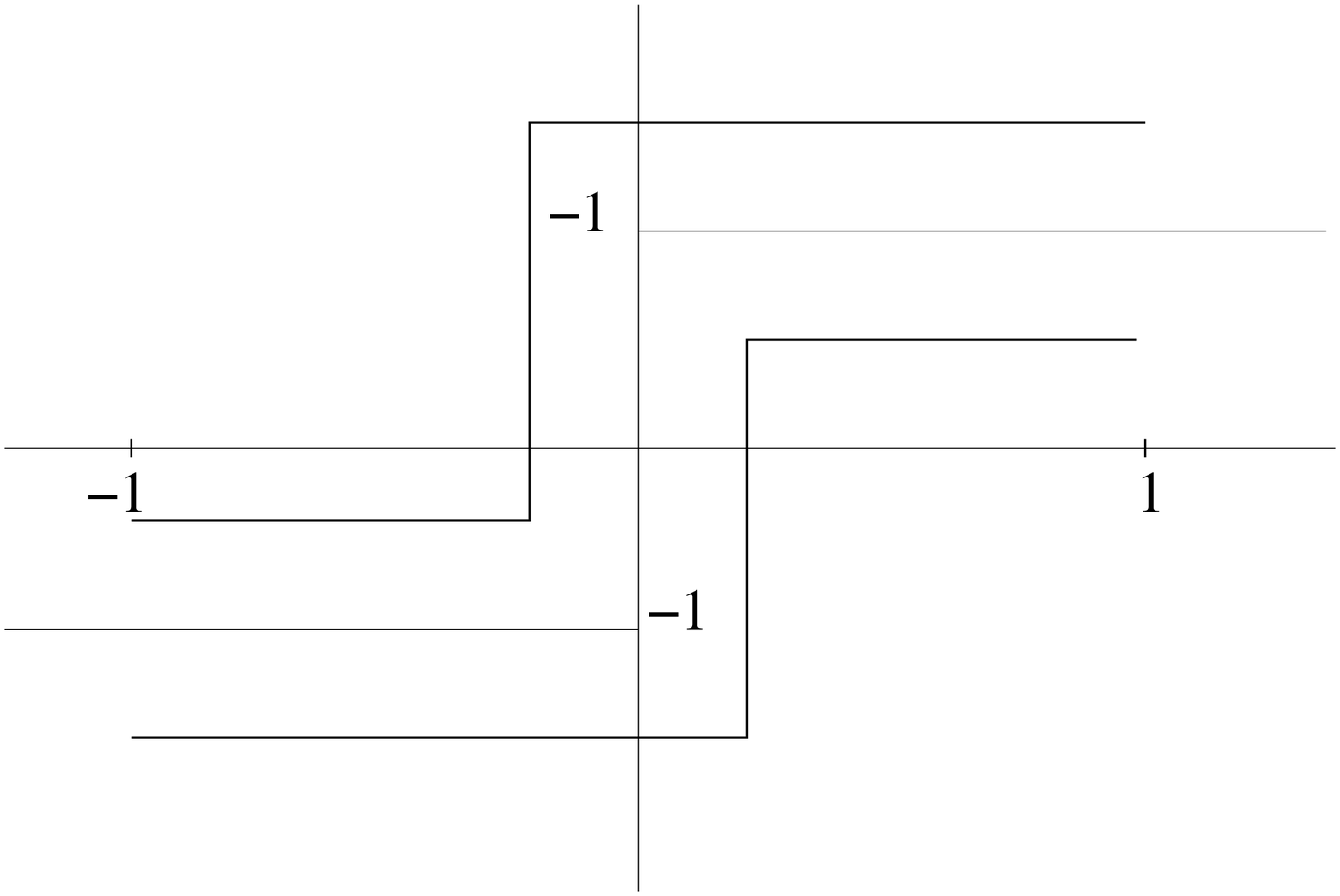}}
Fig. 2. L\'evy's neighborhood of the function $\sgn$.
\end{center}
It is easy to see that if $a=L_m(a)$
then our polynomial $p_m$ from Theorem~1 is the unique polynomial of degree
$2m+1$ which minimizes $\L(\sgn,p)$. We have
\vspace{.1in}

\noindent
{\bf Theorem 5}  
$$\lim_{m\to\infty}\frac{m}{\log m}\L(\sgn,p_m)=\frac{1}{2}.$$
\vspace{.1in}

{\em Remarks}. One could also use the Hausdorff distance between the
completed graphs.
In the case of $\sgn(x)$, the Hausdorff distance will differ
from $\L(\sgn,p_m)$ by a factor of~$\sqrt{2}$. Approximation
of functions with respect to Hausdorff distance between their completed 
graphs was much studied by Sendov \cite{Sendov} and his followers.
An {\em arbitrary} bounded function on $[0,1]$ can be approximated
in this sense by polynomials of degree $n$ with error $O((\log n)/n)$
\cite{Sendov2}.
\vspace{.1in}

%
\vspace{.1in}

{\em Proof of Theorems 3 and 4}. Let $\phi$ be either $\phi_B$ or $\phi_{B,m}$.
By inspection of the boundary correspondence, we conclude that
$f=1-L\cos\phi$ is real on the positive ray and pure imaginary
on the positive imaginary ray. So $f$ extends to an entire function by
two reflections. The extended function evidently satisfies
$$\overline{f(\overline{z})}=
f(z)\quad\mbox{and}\quad-\overline{f(-\overline{z})}=f(z),$$
so we conclude that $f$ is odd. In the case of Theorem~1, the region
$\Omega_{B,m}$ is close to the strip $\{ z:\Re z\in (\pi/2,\pi(m+1))\}$
as $\Im z\to\infty$, so $\phi_{B,m}\sim i(2m+1)\log z,\; z\to\infty$,
so $f$ is a polynomial of degree $2m+1$. In Theorem~2, $\phi(z)\sim z$,
so $f$ has exponential type one.
In both theorems, differentiation shows that the only critical points
of $f$ in the closed right half-plane are preimages of the critical
points of the cosine under $\phi$. So the graph of $f$ has the required
shape. That a polynomial with such graph is the unique extremal
for Theorem~1 follows from the general
theorem of Chebyshev on the uniform approximation of
continuous functions \cite[Ch. II]{Akha}.

The proof that the entire function $f$ we just constructed
is the unique extremal for Theorem~2 might not be so well-known,
so we include this proof which we learned from B. Ya. Levin
(compare \cite{Ere,SY}).

All critical points of our entire function $f$ are real.
Let $x_1<x_2<\ldots$ be the sequence of positive critical points of $f$,
Then we have $x_1>a$, and 
\begin{equation}
\label{crit}
f(x_k)=1+(-1)^{k-1}L,\quad\mbox{and}\quad f(A)=1-L.
\end{equation}
Let $g$ be another real entire function of exponential type $1$
such that 
\begin{equation}
\label{sup}
\sup|g(x)-1|\leq L\quad\mbox{for}\quad x\geq A.
\end{equation}
We may assume without loss of generality that $g$ is odd
(otherwise replace it by $(g(x)-g(-x))/2$ which also
satisfies (\ref{sup})).
Equations (\ref{crit}) and (\ref{sup})
imply that the graph of $g$ intersects the graph of $f$ on
every interval $[x_k,x_{k+1}],\, k\geq 1$. More precisely,
there is a sequence of zeros $y_k$ of $f-g$ (where multiple zeros
are repeated according to their multiplicity), which is interlacent
with $x_k$, that is 
$$x_1\leq y_1\leq x_2\leq y_2\leq\ldots,$$
and {\em in addition} to those $y_j$, $f-g$ has at least one zero
in $(0,x_1]$. 
By the well-known theorem \cite[VII, Thm. 1]{Levin}, it follows that
the meromorphic function
$$F(z)=\frac{1}{z}\prod_{k\in\Z\backslash\{0\}}\frac{1-z/x_k}{1-z/y_k}$$
has imaginary part of constant sign
in the upper half-plane, and of opposite sign in the lower half-plane.
This implies that
\begin{equation}
\label{O}
F(re^{i\theta})=O(r),
\end{equation}
when $r\to+\infty$, uniformly with respect to $\theta$ 
for $\epsilon\leq\theta\leq\pi-\epsilon$, for every $\epsilon>0$.
Similar estimate holds in the lower half-plane.

As $(f-g)(y_k)=0$, we have 
\begin{equation}
\label{77}
\frac{f-g}{f'}=P/F,
\end{equation}
where $P$ is an entire function of exponential type.

It is easy to see that the left hand side of (\ref{77}) is bounded
for $|\Im z|\geq 1$. Indeed,
Phragm\'en and Lindel\"of 
give $|f(z)-g(z)|\leq C_1\exp|\Im z|$,
while $f'$ has only real zeros and approaches
$L\cos(x\pm\alpha))$ as $x\to\pm\infty$, where $\alpha$ is
some real constant. It follows that $|f'(z)|\geq C_2 \exp|\Im z|$
for $|\Im z|>1,$ so $(f-g)/f'$ is bounded for $|\Im z|>1.$

So we conclude from (\ref{O}) that $P(z)=O(|z|)$
and this contradicts the fact that $P$ has at least two zeros,
unless $P=f-g=0$. This completes the proof.
\vspace{.2in}

{\em Proof of Theorem 1}. 
\nopagebreak
We recall that $B_m=B_m(a),\,\Omega_m=\Omega_{B_m,m}$
and the conformal maps of the first quadrant $Q$ onto $\Omega_m$
were defined before Theorem~4. We are going to prove (\ref{A}) first,
which is the same as 
\begin{equation}
\label{11}
B_m=\left( m+\frac{1}{2}\right)\log\frac{1+a}{1-a}+\frac{1}{2}\log m+
\frac{1}{2}\log\frac{2a}{1-a^2}+c+o(1),
\end{equation}
as $m\to\infty$ and $a$ is fixed. Here $c$ is an absolute constant.

We need some auxiliary conformal maps. Let $\psi:Q\to Q$ be
defined by
$$\psi(z)=\frac{A_m}{a}\sqrt{\frac{z^2-a^2}{1-z^2}},$$
where
\begin{equation}
\label{A_m}
A_m=\left( m+\frac{1}{2}\right)\log\frac{1+a}{1-a}>0
\end{equation}
The reason for such choice of $A_m$ will be seen later.
Then $\psi$ gives the following boundary points correspondence:
$$\psi:(0,a,1,\infty)\mapsto(iA_m,0,\infty,iA_m/a).$$
The function $\Phi_m(z)=i\phi_m\circ\psi^{-1}(-iz)$
maps the {\em second} quadrant onto $i\Omega_m$, and sends the positive
imaginary axis onto the interval $\ell=(0,i\pi(m+1))$.
By reflection, we extend $\Phi_m$ to a map from the upper half-plane
onto the region $i(\Omega_m\cup\overline{\Omega_m}\cup\ell)$.
This extended map $\Phi_m$ gives the following boundary correspondence:
$$\Phi_m:(-C_m,-A_m,0,A_m,C_m)\mapsto(-\infty,-B_m,0,B_m,+\infty),$$
where we set $C_m=A_m/a$.

Now we introduce the map $h_m(z)=\Phi_m(z+A_m)-B_m.$
Our first goal is to show that the sequence $h_m$ tends to a limit,
and to describe this limit. The boundary correspondence under $h_m$
is this:
$$h_m:(-C_m-A_m,-2A_m,-A_m,0,C_m-A_m)\mapsto(-\infty,-2B_m,-B_m,0,+\infty).$$
We represent $h_m$ as the Schwarz integral of its imaginary part:
\begin{equation}
\label{schwarz}
h_m(z)=\left(m+\frac{1}{2}\right)\log\frac{1+za/((1+a)A_m)}{1-za/((1-a)A_m)}+
\frac{1}{\pi}\int_{-\infty}^\infty\left(\frac{1}{t-z}-\frac{1}{t}\right)v_m(t)dt,
\end{equation}
where
$$v_m(t)=\left\{\begin{array}{ll}\Im h_m(t),& t\in[-C_m-A_m,\, C_m-A_m],\\
                                 \pi/2,     & t\notin[
-C_m-A_m,\, C_m-A_m].
\end{array}\right.$$
Our choice of $A_m$ in (\ref{A_m}) implies that the first summand in
the right hand side of (\ref{schwarz}) has a limit
$$\lim_{m\to\infty}\left( m+\frac{1}{2}\right)
\log\frac{1+za/((1+a)A_m)}{1-za/((1-a)A_m)}=\sigma z,$$
where
\begin{equation}\label{sigma}
\sigma=\frac{2a}{(1-a^2)\log((1+a)/(1-a))}.
\end{equation}
It is easy to see that the integral in (\ref{schwarz}) converges to
a bounded function of the form
$$\frac{1}{\pi}\int\left(\frac{1}{t-z}-\frac{1}{t}\right)\rho_\sigma(t)dt,$$
where $\rho_\sigma$ is a bounded positive function which we will describe
shortly.
The image of $h_m$ has a limit $\Omega^*$ in the sense of
Caratheodory; $\Omega^*$ is the region in the upper half-plane
above the graph of the function $\arccos e^{-x},\, x\geq 0.$
So $h_m\to H_\sigma$, where $H_\sigma$ is the conformal map
of the upper half-plane onto $\Omega^*$, $H_\sigma(0)=0$
and $H_\sigma(z)\sim\sigma z$ as $z\to\infty.$ 

Thus we have a
Schwarz representation
$$H_\sigma(z)=\sigma z+\frac{1}{\pi}\int\left(\frac{1}{t-z}-\frac{1}{t}\right)
\rho_\sigma(t)dt,$$
and $\rho_\sigma(t)=\Im H_\sigma(t)$ for real $t$.
Our next goal is to study asymptotics of $B_m=-h_m(-A_m)$.
We use the following comparison function
$$
g_m(z)=\left( m+\frac{1}{2}\right)\log\frac{1+za/((1+a)A_m)}{
1-za/((1-a)A_m)}
+\frac{1}{2}\int_0^\infty\left(\frac{1}{t-z}-\frac{1}{t}\right)
\chi_m(t)dt,$$
where $\chi_m$ is the characteristic function
of the set
$(-\infty,-2(A_m+1)]\cup[2,+\infty)$. 
We have
$$\lim_{m\to\infty}(h_m(-A_m)-g_m(-A_m))=-\frac{1}{\pi}
\int_0^\infty\left( \rho_\sigma(t)-
\frac{\pi}{2}\chi_{[2,\infty)}(t)\right)\frac{dt}{t},$$
and using (\ref{A_m}),
$$g_m(-A_m)=-A_m-\frac{1}{2}\log(A_m+1).$$
Combining these two equations we obtain
$$B_m=-h_m(-A_m)=A_m+\frac{1}{2}\log A_m+\frac{1}{\pi}\int_0^\infty\left(
\rho_\sigma(t)-\frac{\pi}{2}\chi_{[2,\infty)}(t)\right)\frac{dt}{t}+o(1).$$
using the evident transformation law $H_\sigma(\lambda z)=H_{\sigma\lambda}(z)$ 
we obtain $\rho_\sigma(\lambda t)=\rho_{\sigma\lambda}(t)$, and therefore,
$$B_m=A_m+\frac{1}{2}\log A_m+\frac{1}{2}\log\sigma+\frac{1}{\pi}
\int_0^\infty\left(\rho_1(t)-\frac{\pi}{2}\chi_{[2,\infty)}(t)\right)
\frac{dt}{t}+o(1).$$
Substituting the values of $A_m$ and $\sigma$ from (\ref{A_m}) and
(\ref{sigma}) we obtain (\ref{11}) with
\begin{equation}\label{c}
c=\frac{1}{\pi}\int_0^\infty\left(\rho_1(t)-\frac{\pi}{2}
\chi_{[2,\infty)}\right)\frac{dt}{t}.
\end{equation}

The numerical value 
$$c=(1/2)\log(2\pi)$$
is
obtained from comparison with Bernstein's result (\ref{bern}).

Indeed, in view of (\ref{our}), (\ref{BL}) and (\ref{11})
we have
\begin{eqnarray*}
L_m(a)&=&\inf_{\deg q=m}\,\sup_{[a^2,1]}\sqrt{x}|q(x)-1/\sqrt{x}|\sim
2e^{-B_m}\\ \\
&=&
2(e^{-c}+o(1))\frac{1-a}{\sqrt{2a}}\left(\frac{1-a}{1+a}\right)^m
\frac{1}{\sqrt{m}}.
\end{eqnarray*}
Comparing this expression for $L_m(a)$ with the expression
(\ref{bern}) when $a\to 1$, we obtain (\ref{11}) with $c=(1/2)\log(2\pi)$. 

\vspace{.2in}

{\em Proof of Theorem 2}.
We have to prove that
$B=A+(1/2)\log A+c_0+o(1)$ as $A\to\infty$, where $c_0=(1/2)\log(\pi/2).$

 Let $f_1(z)=\sqrt{z^2-A^2}$ be the conformal map
of the first quadrant $Q$ onto itself, sending $A$ to $0$,
and $f_1(z)\sim z$ as $z\to\infty$. Then $f_1(0)=iA.$
Let $f_2$ be the conformal map of $Q$ onto $\Omega_B$, $f_2(0)=0$
and $f_2(z)\sim z,$ as $z\to\infty$. We extend $f_2$ by symmetry,
reflecting both domains in the positive ray. So from now on
$f_2$ is defined in the right half-plane. The condition that $f_2(iA)=iB$
defines the number $B$ uniquely. It is easy to see that
$$B\sim A$$
as $A\to\infty$.

Now put $h(z)=if_2(-iz)$ for convenience.
This $h$ maps the upper half-plane onto a subregion of the upper
halfplane. The boundary of this subregion is asymptotic to
the line $\{\Im z=\pi/2\}$ as $|\Re z|\to\infty $.

We have $B=h(A)$,
and we wish to find asymptotics of $h(A)$ as $A\to\infty$.

To do this, we use the following comparison function
$$g(z)=z+\int_{A+2}^\infty\frac{z}{t^2-z^2}dt.$$
The integral in the right hand side is the 
Schwarz formula for an analytic function in the upper half-plane
whose imaginary part equals $$(\pi/2)\chi_{(-\infty,-A-2]\cup
[A+2,\infty)}.$$ 

Our function $h$ has a similar representation in terms of its
imaginary part on the real line. Subtracting these two
representations, we obtain
$$g(A)-h(A)=\frac{2A}{\pi}\int_0^\infty\frac{v(t+A)}{2At+t^2}dt,$$
where $v(t)=\Im (g(t)-h(t))$.
Now we claim that $v(t+A)$ tends
to a limit $v_0$, as $A\to+\infty$.
This limit is
$(\pi/2)\chi_{[2,\infty)}(t)-\Im H(t)$,
where 
and $H(z)=\lim_{A\to+\infty}( h(z+A)-B(A))$
is the conformal map of the upper half-plane onto
the region in the upper half-plane above the graph
$y=\arccos e^{-x},\, x\geq 0$ and $y=0,\, x\leq 0$.
Notice that $H=H_1$, where $H_1$ was defined in the
proof of Theorem 1. 

Thus $g(A)-h(A)\to\const$ where the constant is given by
\begin{equation}
\label{C_0}
\frac{1}{\pi}
\int_0^\infty\frac{v_0(t)}{t}dt=\frac{1}{\pi}\int_0^\infty
(\frac{\pi}{2}\chi_{[2,\infty)}(t)-\Im H(t))\frac{dt}{t}.
\end{equation}
The integral is convergent because $v_0(t)=O(\sqrt{t}),\; t\to 0$
and $v_0(t)=O(e^{-t}),\; t\to\infty$.
It remains to find the asymptotic behavior
of $g(A)$ as $A\to\infty$.
We have
$$g(A)=A+\int_2^\infty\frac{A}{2tA+t^2}dt=
A+\frac{1}{2}\log(A+1).$$
Combining these results we obtain
\begin{eqnarray*}
B=h(A)&=&
A+\frac{1}{2}\log A+\frac{1}{\pi}\int_0^\infty\left(\Im H(t)-
\frac{\pi}{2}\chi_{[2,\infty)}(t)\right)\frac{dt}{t}+o(1)\\ \\
&=& A+\frac{1}{2}\log A+c+o(1),
\end{eqnarray*}
where $c$ is the constant from (\ref{c}).
This proves Theorem 2. 
\vspace{.2in}

{\em Proof of Theorem 5}. We choose 
\begin{equation}
\label{tri}
B=\log m-\log\log m+\log 4,
\end{equation}
so that $L\sim(\log m)/(2m)$ in view of (\ref{BL}).
Consider the conformal map of $\Omega_{B,m}$ onto the half-strip
$\Pi=\{ z:\Re z\in(0,(m+1)\pi)\}$ by a function $f_1$ such that
$f_1(0)=0,\, f_1((m+1)\pi)=(m+1)\pi$ and $f_1(\infty)=\infty$.
Let $B'=f_1(B)$.
The function $f_1(z/B)$ tends to the identity as $m\to\infty$, so
$$B'\sim B,\quad B\to\infty.$$
Let $f_2$ be the (elementary) conformal map of $\Pi$ onto the first quadrant
normalized by $f_2(B')=0,\; f_2(\pi(m+1))=1$ and $f_2(\infty)=\infty$. Then
$$a:=f_2(0)=\frac{1-\exp(B'/(m+1))}{1+\exp(B'/(m+1))}=
B'/2(m+1)\sim B/2m\sim (\log m)/(2m).$$

The authors thank Doron Lubinsky, Misha Sodin and Andrei Gabrielov
for their help and useful comments.

{\em A. E.: eremenko@math.purdue.edu
\nopagebreak

Department of Mathematics, Purdue University,
\nopagebreak

West Lafayette, IN 47907-2067
\nopagebreak

U. S. A.

\vspace{.1in}

P. Yu.: yuditski@macs.biu.ac.il

Department of Mathematics, Bar Ilan University,

52900 Ramat Gan,

Israel}
\end{document}